\documentclass{amsart}
\usepackage{amssymb}
\usepackage[square]{natbib}
\usepackage{epsfig}
\newtheorem{mainthm}{Theorem}[]\newtheorem{maincor}[mainthm]{Corollary}
\newtheorem{thm}{Theorem}[section]\newtheorem*{conj}{Conjecture}
\newtheorem*{prob}{Problem}

\newtheorem*{thm*}{Theorem}   
 
\newtheorem{lem}[thm]{Lemma}

\newtheorem{prop}[thm]{Proposition}
\newtheorem{sublem}[thm]{Sublemma}
\theoremstyle{definition} 
\newtheorem{step}{Step}
\newtheorem{mainrem}[mainthm]{Remark}
\newtheorem*{rems*}{Remarks} 
\theoremstyle{remark}
\newtheorem{rem}[thm]{Remark}

\newcommand{\fF}{\mathbb{F}}
     
 \newcommand{\fD}{\mathbb{D}}

\newcommand{\R}{\mathbb{R}}

\newcommand{\folL}{\mathcal{L}}
\newcommand{\folF}{\mathcal{F}}

\newcommand{\tc}{\tilde{c}}

\newcommand{\G}{\mathsf{G}}\newcommand{\bgH}{\bar{\mathsf{H}}}

\newcommand{\gH}{\mathsf{H}}

\newcommand{\dc}{\dot{c}}\newcommand{\gD}{\mathsf{D}}
\DeclareMathOperator{\Par}{Par}

\DeclareMathOperator{\dis}{dis}

\DeclareMathOperator{\Or}{O}

\DeclareMathOperator{\Diff}{Diff}

\newcommand{\gK}{\mathsf{K}}

\DeclareMathOperator{\pr}{pr}

\DeclareMathOperator{\spann}{span}
 
 \newcommand{\bbJ}{\Lambda} 
 \newcommand{\bbV}{\Upsilon}
     
\newcommand{\ml}{\langle}                     
\newcommand{\mr}{\rangle}                     

\renewcommand{\fF}{\mathcal{F}}
\newcommand{\fL}{\mathcal{L}}

\newcommand{\eps}{\varepsilon}

\begin{document}

%
\begin{titlepage}\title[Riemannian foliations in nonnegative curvature]
{A duality theorem for Riemannian foliations in nonnegative sectional
 curvature}
\author{Burkhard Wilking}
\thanks{to appear in GAFA}
\end{titlepage}
\maketitle
\setcounter{page}{1}
\setcounter{tocdepth}{0}

Using a new type of Jacobi field estimate 
 we will prove a 
 duality theorem for singular 
Riemannian foliations in complete manifolds of nonnegative sectional 
curvature.  Recall that a transnormal system $\fF$ 
is a subdivision of $M$ 
into $C^\infty$ immersed connected complete submanifolds 
without boundary, called leaves, 
such that geodesics emanating perpendicularly to one 
leaf stay perpendicularly to the leaves. 
If $M$ is complete the leaf  $\fL(p)$ of each point $p\in M$ 
is intrinsically complete as well.
A transnormal system $\fF$ is called a singular Riemannian 
foliation if there are vectorfields $X_i$ ($i\in I$) in $M$ 
such that  $T_p \folL(p)=\spann_{\R}\{X_{i|p}\mid i\in I\}$ for all $p\in M$, 
see [Molino, 1988].
Examples of singular Riemannian foliations are the fiber decomposition
of a Riemannian submersion or the orbit decomposition
 of an isometric group action.

A piecewise smooth curve $c$ is called horizontal 
with respect to a transnormal system 
$\fF$, if $\dc(t)$ is in the normal bundle 
$\nu_{c(t)}\bigl(\folL(c(t))\bigr)$ of the leaf $\fL(c(t))$.
One can define a dual foliation $\folF^\#$ by defining
\[
\fL^\#(p):=\bigl\{q\in M\mid \mbox{ there is a piecewise 
smooth horizontal curve from $p$ to $q$}\bigr\}
\] 
as dual leaf of a point $p\in M$.
We will see that $\fL^\#(p)$ is a smooth immersed submanifold 
of $M$, see section~\ref{sec: general}.
The double dual is not always equal to the original 
foliation. But the triple dual foliation is usually isomorphic to 
the dual foliation.
In general one cannot expect that the dual foliation has
too many reasonable properties. 
We will see that this is different in nonnegative curvature. 
The main results can be interpreted as 
rigidity versions of the following 

\begin{mainthm}\label{thm1}\label{thm: K>0}
 Suppose that $M$ is a complete positively curved manifold 
with a singular Riemannian  foliation $\fF$. Then the dual foliation 
has only one leaf.
\end{mainthm}

In other words one can connect two arbitrary points in $M$ by
 a horizontal curve. It should be noted that the theorem 
is global in nature.
 If one considers a cohomogeneity one action on a sphere 
then of course the horizontal distribution is 
one dimensional in the generic part and hence integrable. 
However, a horizontal curve can run into singular 
orbits and then with different directions 
out of singular orbits. This way one can reach more than 
just  a one dimensional subset and in fact every point on the sphere.

Theorem~\ref{thm1} suggests to introduce a new length metric on 
$M$ by defining the distance of two points as the infimum over the length 
of all horizontal curves connecting these two points. 
The previous example shows that one cannot expect that 
the two metrics induce the same topology, but it would 
be interesting to know whether one can say more 
about the latter metric, other than that $M$ stays connected.

We prove Theorem~\ref{thm1} together with the following rigidity result 
in section~\ref{sec: thm1}.

\begin{mainthm}\label{thm2}\label{thm: micro}
Suppose that $M$ is a complete nonnegatively curved manifold 
with a singular Riemannian  foliation $\fF$. Suppose the leaves of
 the dual foliation are complete. 
Then $\fF^\#$ is a singular Riemannian foliation as well. 
\end{mainthm}

In many cases it is actually possible 
to remove the assumption on the completeness of the dual leaves.

\begin{mainthm}\label{thm3}\label{thm: complete}
 Suppose that $M$ is a complete nonnegatively curved manifold 
with a singular Riemannian  foliation $\fF$. 
Then the dual foliation has intrinsically complete leaves if in addition 
one of the following holds. 
\begin{enumerate}
\item[a)] $\fF$ is given by the orbit decomposition of an 
isometric group action. 
\item[b)] $\fF$ is a non singular foliation and $M$ is compact.
\item[c)] $\fF$ is given by the fibers of the Sharafutdinov retraction.
\end{enumerate}
\end{mainthm}

We recall that
 an open nonnegatively curved manifold $M$ is by the soul theorem 
of Cheeger and Gromoll [1972] diffeomorphic 
to the normal bundle of a compact totally geodesic submanifold $\Sigma$,
 the soul of $M$. Sharafutdinov showed that there is a distance 
nonincreasing retraction $P\colon M\rightarrow \Sigma$. 
By Perelman's [1994] solution of the soul conjecture 
$P$ is a Riemannian submersion of class $C^{1}$. 
Guijarro [2000] improved the regularity to $C^2$. 
Before we prove Theorem~\ref{thm: complete} in section~\ref{sec: complete}
we will use  Theorem~\ref{thm: micro} 
to establish the following regularity result in section~\ref{sec: smooth}.

\begin{maincor}\label{cor: smooth} Let $(M,g)$ be an open nonnegatively curved 
manifold,  $\Sigma$ a soul of $M$. Then 
the Sharafutdinov retraction 
$P\colon M\rightarrow \Sigma$ is of class $C^\infty$.
\end{maincor}

  Cao and Shaw [2005] proposed an independent proof
of Corollary~\ref{cor: smooth}. They showed that the fibers
of the Sharafutdinov retraction admit locally an one dimensional foliation
by geodesics. However, the tangent fields of these geodesics 
 are obtained from the convex exhaustion in the soul construction
and
it is not at all clear that this foliation of geodesics
 is of class $C^\infty$,  even in regions where $P$ is of class $C^\infty$. 
Thus the same holds for the maps constructed with this foliation. 
This causes crucial problems in the proof of Proposition 6
 in [Cao and Shaw, 2005] which one can probably not overcome.

Recall that a map  between metric spaces $\sigma\colon X\rightarrow Y$ 
is a submetry if $\sigma (B_r(p))=B_r(\sigma(p))$ 
for any metric ball $B_r(p)$ in $X$. If both 
$X$ and $Y$ are Riemannian manifolds then $\sigma$ 
is a Riemannian submersion of class $C^{1,1}$ by a result of 
Berestovskii and Guijarro [2000]. If 
 $X$ is a Riemannian manifold and $Y$ is arbitrary, then 
the fibers of 
$\sigma$ give rise to  a generalized singular Riemannian foliation.  
In general the fibers can have boundary and might be only of 
class $C^{1,1}$.

\begin{maincor}\label{cor: product submetry}\label{cor: submetry} 
Let $(M,g)$ be an open nonnegatively curved 
manifold, and let  $\Sigma$ be a soul of $M$. Then there is 
a noncompact Alexandrov space $A$ and a submetry 
\[
\sigma\colon M\rightarrow \Sigma \times A
\] 
where $ \Sigma \times A$ is endowed with the product metric. 
Moreover the fibers of $\sigma$ are compact smooth
submanifolds without boundary.  
\end{maincor}

In particular, any non-contractible open 
manifold of nonnegative sectional curvature
has a nontrivial product as a metric quotient, for the proof 
see section~\ref{sec: product}.

\begin{maincor}\label{cor: flats}
 Let $(M,g)$ be an open nonnegatively curved manifold,
  $\Sigma$ a soul of $M$, and 
$P\colon M\rightarrow \Sigma$ the Sharafutdinov retraction. 
Suppose $x\in T_pM$ is horizontal with respect 
to $P$, and suppose that $v\in T_pM$ is a vertical vector 
perpendicular to the holonomy orbit. Then $v$ and $x$ 
span a totally geodesic flat.
\end{maincor}

We should mention that the family of totally geodesic flats 
in Corollary~\ref{cor: flats} is at least
as big as the family obtained from Perelman's proof of the soul conjecture.
Equality occurs precisely if the normal holonomy group of 
$\Sigma$ acts transitively on the normal sphere.
From a metric point of view this is a somewhat special 
case which is better understood than the general situation.
For example by Walschap [1988] the normal exponential map of the soul 
is a diffeomorphism  and the 
cone at infinity is a ray. 
In general the diffeomorphism between $\nu(\Sigma)$ 
and $M$  is not given by the exponential map. 
However using our results we can show that the diffeomorphism 
can be chosen 
such that it respects the structure of $M$ as a doubly foliated 
space, see section~\ref{sec: linear}. 
\begin{maincor}\label{cor: linearization}
 There is a diffeomorphism $f\colon \nu(\Sigma)\rightarrow M$
satisfying
\begin{enumerate}
\item[a)] $P\circ f=\pi$, where $\pi\colon \nu(\Sigma)\rightarrow \Sigma$ 
denotes the natural projection. 
\item[b)] $f_*$ maps the horizontal geodesics in $\nu(\Sigma)$
onto horizontal geodesics  in $M$, where $\nu(\Sigma)$ is endowed
with the natural connection metric.  
\end{enumerate}
\end{maincor}

 It was shown in
[Wi, 2004] that if a group $\G$ acts isometrically on a positively curved
 manifold with a nontrivial principal isotropy group, 
then the orbit space $M/\G$ has boundary. 
In nonnegative sectional curvature we have
 the following rigidity result (section~\ref{sec: actions}). 

\begin{maincor}\label{cor: actions}
 Let $(M,g)$ be a  nonnegatively curved complete 
manifold, and suppose a Lie group $\G$ 
acts isometrically and effectively 
 on $(M,g)$ with principal isotropy group $\gH\neq 1$. 
If the orbit space $M/\G$ has no boundary, then 
there is a closed subgroup $\gK$ 
with $\gH\triangleleft \gK$,
an invariant metric on $\G/\gK$, and a $\G$ 
equivariant Riemannian submersion 
$\sigma\colon M\rightarrow \G/\gK$ with totally geodesic fibers. 
\end{maincor}

The main new tool used to prove the above results 
is a simple and general observation which may very well be useful in 
different context as well. It allows to give what we call transversal 
Jacobi field estimates.
Let $c\colon I\rightarrow (M,g)$ be a geodesic in 
a Riemannian manifold $(M,g)$, and let
$\bbJ$  be an $(n-1)-$dimensional family of normal
Jacobi fields for which the 
 corresponding Riccati operator is self adjoint. 
Recall that the Riccati operator $L(t)$ is the endomorphism of 
$(\dc(t))^\perp$ defined by $L(t) J(t)=J'(t)$ 
for $J\in \bbJ$.
Suppose we have a vector subspace
  $\bbV\subset\bbJ$. 
Put 
\[
T^v_{c(t)}M:=\{J(t)\mid J\in \bbV\}\oplus \{J'(t)\mid J\in \bbV, J(t)=0\}.
\]
Observe that the second summand vanishes for almost every 
$t$ and that  $T^v_{c(t)}M$ depends smoothly on $t$.
We let $T^{\perp}_{c(t)}M$ denote the orthogonal 
complement of   $T^v_{c(t)}M$, and 
for $v\in T_{c(t)}M$ we define $v^{\perp}$ 
as the orthogonal projection of $v$ to $T^{\perp}_{c(t)}M$. 
If $L$ is non-singular at $t$ we put 
\[
A_t\colon T^v_{c(t)}M\rightarrow T^{\perp}_{c(t)}M, \mbox{ 
 $J(t)\mapsto J'(t)^{\perp}$ for $J\in \bbV$}.
\] 
It is easy to see that $A$ can be extended continuously on $I$. 
For a vector field $X(t)\in T^{\perp}_{c(t)}M$ we 
define $\tfrac{\nabla^{\perp}X}{\partial t}= (X'(t))^\perp$.
The following observation which is proved in section~\ref{sec: transversal}
 is key.

\begin{mainthm}\label{thm: jacobi}\label{thm: transverse estimate}
 Let $J\in \bbJ -   \bbV$ and put $Y(t):=J^{\perp}(t)$.
 Then $Y$ satisfies the following Jacobi equation 
\[
\tfrac{(\nabla^{\perp})^2}{\partial t^2}Y(t) +
 \bigl(R(Y(t),\dc(t))\dc(t)\bigr)^{\perp} + 3 A_tA^*_tY(t)=0.
\]
\end{mainthm}

One should consider $\bigl(R(\cdot,\dc(t))\dc(t)\bigr)^{\perp} + 3 A_tA^*_t$ 
as the modified curvature operator.  
The crucial point in the equation is that the additional  O'Neill type term
$ 3 A_tA^*_t$  is positive 
semidefinite. We will denote the family of all 
vector fields $Y$ obtained from $\bbJ$ and $\bbV$ 
with $\bbJ/\bbV$.

\begin{maincor}\label{cor: decomposition}\label{cor: decomp} 
Consider  an $n-1$-dimensional family $\bbJ$ of 
normal Jacobi fields with a self adjoint Riccati operator along a geodesic $c\colon \R \rightarrow M$
 in a
nonnegatively curved manifold. 
Then 
\[
\bbJ=\spann_{\R}\bigl\{ J\in \bbJ \mid \mbox{ $J(t)=0$ for some $t$}\bigr\}\oplus 
\bigl\{ J\in \bbJ\mid \mbox{ $J$ is parallel}\bigr\}.
\]
\end{maincor}

It should be understood that this does 
not follow from the usual Rauch or Riccati comparison for the family $\bbJ$
 since 
 this fails after the first conjugate point. 
Instead one considers 
\[
\bbV:=\{ J\in \bbJ \mid \mbox{ $J(t)=0$ for some $t$}\}.
\]
 Then  for any $J\in \bbJ -  \bbV$ and any $t\in \R$ the vector 
$J(t)$ is transversal to  
\[
T^v_{c(t)}M:=\{J(t)\mid J\in \bbV\}\oplus \{J'(t)\mid J\in \bbV, J(t)=0\}.
\]
By Theorem~\ref{thm: transverse estimate}
 the family  $\bbJ/\bbV$ 
satisfies again a Jacobi equation with nonnegative 
curvature operator, and as explained for this family the selfdual 
Riccati operator is non-singular everywhere. 
By the usual Riccati comparison (see for example 
[Eschen\-burg and Heintze, 1990])
the Riccati operator 
of the family $  \bbJ/\bbV$ vanishes 
and $\bbJ/\bbV$ consists of parallel Jacobi fields.
Clearly Corollary~\ref{cor: decomp} follows. 
We conclude the introduction with a few open problems.

\begin{prob}[Berestovskii and Guijarro] Let $\sigma\colon M\rightarrow B$ 
a submetry between complete nonnegatively curved manifolds. 
Is $\sigma$ of class $C^\infty$?
\end{prob}

If one assumes in addition that  $M$ is compact and that 
$\sigma$ is of class $C^\infty$
on some open subset $U\subset M$, then it is 
conceivable that one can modify the proof of Corollary~\ref{cor: smooth} 
to show that $\sigma$ is smooth.
Similarly Theorem~\ref{thm: complete}
 might be viewed as support for the following

\begin{conj} Suppose $\folF$ is a singular Riemannian foliation of
a nonnegatively curved complete manifold $M$. 
Then the dual foliation has complete leaves.
\end{conj}

In singular spaces one can define metric foliations 
as subdivisions into connected subsets which 
are locally given by the fibers of a submetry. One can still
define horizontal curves in this setting and it is natural 
to ask

\begin{prob} Suppose $X$ is an Alexandrov 
space of nonnegative curvature and suppose 
that $\folF$ is a metric foliation. Is the dual 
foliation also a metric foliation?
\end{prob}

The idea for this paper came when the author thought
about a problem posed by  V.~Kapovitch 
proposing that collapse of manifolds with lower curvature bound
should in a suitable sense occur along the fibers of a submetry. 
If $M$ is an open manifold of nonnegative sectional curvature, 
then the cone at infinity $C(M)$ of $M$  
is isometric to the cone at infinity of $A$
where $A$ is the Alexandrov space from 
Corollary~\ref{cor: product submetry}.
By combining with Perelman's stability theorem, 
 if $\dim(C(M))=\dim(A)$,
then the collapse of $M$ to $C(M)$ indeed occurs along the fibers 
 the submetry $\pr_2\circ\sigma\colon M\rightarrow A$
 from Corollary~\ref{cor: product submetry}.

\begin{mainrem}
Corollary~\ref{cor: decomp} gives 
also obstructions for invariant 
positively curved metrics on cohomogeneity one 
manifolds. In fact, if $c(t)$ is a normal geodesic in a positively curved 
cohomogeneity one $\G$--manifold, then the Killing fields of the 
action give rise to an $(n-1)$--dimensional family of 
Jacobi fields along $c$ with a self adjoint Riccati operator.
Applying Corollary~\ref{cor: decomp} gives that the Lie algebras 
of the isotropy groups along $c$ generate the Lie algebra of $\G$
 as a vectorspace. For more details we refer the reader to 
[GWZ].
\end{mainrem}

I would like to thank one of the referees for useful comments.

\section{The Transversal Jacobi Field Equation.}\label{sec: transversal}

In this section we prove Theorem~\ref{thm: transverse estimate}. 
 It suffices to prove the equality for a generic $t_0$, 
i.e. we may assume that the Riccati operator is non-singular at $t_0$ 
or equivalently $(\dc(t_0))^\perp=\{ J(t_0) \mid J\in \Lambda\}$.
 Since we can add a Jacobi field 
of $\bbV$ to $J$ without changing $Y$, we can without loss of generality 
assume that $J(t_0)\in T^\perp_{c(t_0)}M$.
Let $X_1(t),\ldots,X_d(t)\in T^{\perp}_{c(t)}M$
 be  orthonormal vector fields with $\tfrac{\nabla^{\perp}}{\partial t}X_i=0$.
We may assume $J(t_0)=X_1(t_0)$. 
We claim
\begin{equation}\label{eq: one}
(J'(t_0))^v= A^* J(t_0)
\mbox{ and } X'_i(t)=- A^* X_i(t).
\end{equation}
To prove these equations let $V(t)$ denote a Jacobi field in $\bbV$. 
Then 
\[
\ml J'(t_0), V(t_0)\mr = \ml J(t_0), V'(t_0)\mr=\ml J(t), A V(t)\mr,
\]
where we used that the Riccati operator of the family $\bbJ$ 
is self adjoint. 
The second equation of ~\eqref{eq: one} follows from 
\begin{eqnarray*}
0= \tfrac{d}{dt}\ml X_i(t), V(t)\mr&=&
 \ml X_i'(t), V(t)\mr +\ml X_i(t), V'(t)\mr\\
&=& \ml X_i'(t), V(t)\mr +\ml X_i(t), AV(t)\mr.
\end{eqnarray*}

Thus we can finish the proof of Theorem~\ref{thm: transverse estimate} as follows.

\begin{eqnarray*}
\ml \tfrac{(\nabla^{\perp})^2}{\partial t^2}J^\perp, X_k(t_0)\mr\!&+&\! 
\bigl\ml R(J(t_0),\dc(t_0))\dc(t_0), X_k(t_0)\bigr\mr =\\
&=& 
 \tfrac{d^2}{dt^2}_{t=t_0}\ml J, X_k\mr-\ml J''(t_0),X_k(t_0)\mr\\
&=& 2\ml J'(t_0),X_k'(t_0)\mr + \ml J(t_0),X_k''(t_0)\mr\\
&=&  2\ml J'(t_0), X_k'(t_0)\mr + \ml X_1(t_0),X_k''(t_0)\mr\\
&=&  2\ml J'(t_0), X_k'(t_0)\mr - \ml X_1'(t_0),X_k'(t_0)\mr+
\tfrac{d}{dt}_{|t=t_0} \ml X_1(t),X_k'(t)\mr\\
&=& -3\ml AA^* J(t_0), X_k(t_0)\mr, \\
\end{eqnarray*}
where we used $ \ml X_1(t),X_k'(t)\mr\equiv 0$ 
and equation~\eqref{eq: one} for the last equality.
Clearly the theorem follows.

\section{Some General Remarks on Dual Foliations}\label{sec: general}

\begin{prop}\label{prop: general dual}\label{prop: general}
 Let $(M,g)$ be a Riemannian manifold with 
a transnormal system $\folF$, and let $\folF^\#$ denote 
the dual foliation.
 There is a family of $C^\infty$ 
vector fields $(X_i)_{i\in I}$ with compact supports such that 
 any dual leaf $\folL^\#$ is a $C^\infty$ immersed submanifold 
with $T_p \folL^\#=\spann_{\R}\{X_{i|p}\mid i\in I\}$.
\end{prop}

\begin{rem}\label{rem: general}\begin{enumerate}
 \item[a)] Even if the ambient manifold and the leaves of 
$\folF$ are complete it is in general not clear that the dual 
leaves are complete. The dual foliation could for example have open 
leaves.
\item[b)] The proof below also shows that one can connect two 
points of one dual leaf by a piecewise horizontal geodesic.
\end{enumerate}
\end{rem}

\begin{proof}[Proof of Proposition~\ref{prop: general dual}.]
Let $\folL$ be a leaf of $\folF$. By assumption 
 geodesics emanating perpendicularly to $\folL$ 
stay horizontal. 
We consider 
all $C^\infty$ vector fields $X$ 
in $M$ which can be obtained as follows:
there is a relatively compact open  subset $K$ of the normal bundle 
of $\folL$ such that $\exp_{K}$ is a diffeomorphism 
onto its image, the image contains the support of $X$
  and $\exp^{-1}_*(X)$ is a vector field 
tangential to the fiber direction, that 
is $\exp^{-1}_*(X)$ is in the kernel of $\pi_*$ 
where $\pi\colon \nu(\folL)\rightarrow \folL$ 
denotes the natural projection. 
Since the set $\exp(\nu_p(\folL))$ is contained in a dual 
leaf for each $p\in\folL$, it follows that 
 the flow lines of such a vector field 
 stay in a dual leaf. In fact two points on a flow line 
can be connected by a piecewise horizontal 
geodesic.

We let $C$ denote the collection of 
all these vector fields where 
$\folL$ runs as well. Let $\gD$ 
denote the diffeomorphism group 
generated by the flows of all vector fields in $C$.
Finally put 
\[
C_2:=\bigl\{\phi_*X_{\circ \phi^{-1}}\mid X\in C , \phi\in\gD\bigr\}.
\] 
Since $\phi\in \gD$ maps the integral 
curves of $X\in C$ to the integral curves of  
$\phi_*X_{\circ \phi^{-1}}$, the group $\gD$ is 
 also the group generated by the flows of the vector fields 
in $C_2$. 
By construction the singular distribution 
spanned by $C_2$ has constant dimension along orbits of $\gD$. 
Using the description of the 
Lie bracket as a Lie derivative, we see that
 Lie brackets of vectorfields in $C_2$ are tangential 
to the distribution. As in the proof of  Frobenius' theorem 
we see that the orbits 
of $\gD$ are smooth submanifolds 
whose tangent space at each point is spanned by vector fields 
in $C_2$. Since these tangent spaces contain 
all vectors which are horizontal with respect to $\folF$, it is clear that 
the orbits of $\gD$ coincide with the dual leaves.
\end{proof}

\section{The Dual Foliation in Nonnegative Curvature.}\label{sec: micro}\label{sec: thm1}
In this section we prove
Theorem~\ref{thm: K>0} and Theorem~\ref{thm: micro}
simultaneously. Let $\fF$ 
be a singular Riemannian foliation of a complete nonnegatively curved manifold.
Suppose the dual foliation has more than one 
leaf.
 Then there is a dual leaf $\folL^\#$ which 
is not open. By Remark~\ref{rem: general} any two points in $\folL^\#$
can be connected by a piecewise horizontal geodesic.

We first plan to show that for any $\folF$-horizontal geodesic $c\colon \R\rightarrow \folL^\#$
the normal space of $\folL^\#$ along $c$ is spanned by parallel Jacobi fields.

We consider the family $\Lambda$ of normal Jacobi fields 
along $c$ that correspond to variations of $c$ by 
geodesics emanating perpendicularly to $\fL(c(0))$ at time $0$.
Clearly, the Riccati operator corresponding to $\Lambda$ 
is self adjoint. 

Consider a Jacobi field  $J\in \Lambda$
 with $J(t_0)=0$ for some $t_0$.
We want to prove
$J'(t_0)\in \nu(\fL(c(t_0))$.
By assumption $J$ is the  variational vectorfield 
 a variation $c_s(t)$ of $c$ by horizontal geodesics. 
Let $Y_i$ ($i\in I$) be a family of vectorfields satisfying 
$\spann_{\R}\{Y_{i|p}\mid i\in I\}=T_p\folL(p)$ for all $p$.  
Since $\dc_s(t_0)$ is perpendicular to $Y_{i|c_s(t_0)}$
and $\tfrac{d}{d s}_{|s=0} c_s(t_0)=J(t_0)=0$,
we get  
\[
0=\tfrac{\partial}{\partial s}_{|s=0}
\ml Y_{i|c_s(t_0)},\dc_s(t_0)\mr= \ml Y_{i|c(t_0)},J'(t_0)\mr
\] 
for all $i$. Thus $J'(t_0)\in \nu(\fL(c(t_0))$.
This shows that $J$ can be written also as a variation 
of $c$ by horizontal geodesics with a fixed value $c(t_0)$ 
at time $t_0$. 
Therefore $J(t)$ is tangential to the dual 
leaf $\folL^\#$ for all $t$.
Hence the vectorfields in 
\[
\Upsilon:=\spann_{\R}\{J\in \Lambda\mid J(t)=0 \mbox{ for some t }\}
\]
are everywhere tangential to the dual leaf $\folL^\#$.

We deduce from Corollary~\ref{cor: decomp} that  $\Lambda$ contains 
a nontrivial subfamily of parallel Jacobi fields
and this completes the proof of Theorem~\ref{thm: K>0}. 

We proceed with the proof of Theorem~\ref{thm: micro}. 
Each Jacobi field $J\in \Lambda$ with $J(0)\in T\folL^\#$ 
is everywhere tangential to $\folL^\#$. Therefore, the 
subspace
$V\subset \Lambda$ of normal Jacobi fields which are everywhere 
tangential to $\folL^\#$ has the maximal possible dimension 
$\dim(\folL^\#)-1$. 
Notice that the decomposition of Corollary~\ref{cor: decomposition} 
necessarily defines two pointwise
 orthogonal families of Jacobi fields. Because of 
$\Upsilon \subset V$ we see that the normal bundle 
of $\folL^\#$ along $c$ is spanned by parallel Jacobi fields. 
The rest of the proof is divided into three steps.

\begin{step}\label{step1} Let $\folL^\#_0$ be a dual leaf
of maximal dimension. Then 
$\folF^{\#}$ induces a (non-singular) Riemannian 
foliation in the $r$-tube $B_r(\folL^\#_0)$  around 
$\folL^\#_0$ for a suitable small $r>0$. 
\end{step}

 By Proposition~\ref{prop: general dual} there 
is an open neighborhood $U$ of $\folL^\#_0$ 
such that $\folF^\#$ induces an actual foliation on $U$. 
We may assume that $U$ decomposes into dual leaves. 
In fact otherwise we can replace $U$ by 
$\bigcup_{\phi \in\gD} \phi(U)$, where $\gD$ 
denotes the group of diffeomorphisms defined in the proof 
of Proposition~\ref{prop: general dual}. 
Let $\folL^\#\subset U$ be a dual leaf, and let 
$N\subset \nu(\folL^\#)$ be an induced leaf of the normal 
bundle of $\folL^\#$. The natural projection $N\rightarrow \folL^\#$ 
is a covering and along $\folF$-horizontal geodesics 
in $\folL^\#$, the submanifold $N$ develops by parallel Jacobi fields. 
Since any two points in a dual leaf
can be connected by a piecewise horizontal geodesic, it follows that 
$N$ consists of vectors of the same length. 
Clearly this shows 
 that $\folF^\#$ is a Riemannian foliation in $U$. 

Choose a point $p\in \folL^\#_0$ and a number $r>0$ 
such that $B_r(p)\subset U$. 
Since the Riemannian foliation $\folF_{|U}$ decomposes 
into intrinsically complete dual leaves, it follows that 
$B_r(\folL^\#_0)\subset U$. This completes the proof of Step~\ref{step1}.

\begin{step}\label{step2} Let $ \folL^\#_0$ be a dual leaf
of maximal dimension, and let 
$N\subset \nu(\folL^\#_0)$ be an induced leaf of the normal 
bundle of $\folL^\#_0$.  There is a unique maximal $s_0\in (0,\infty]$ 
such that $\exp(sN)$ is a dual leaf of maximal 
dimension for all $0<s<s_0$.
If $s_0<\infty$, then $\exp(s_0N)$ is a dual leaf 
whose dimension is not maximal. Furthermore, 
the map $N\rightarrow \exp(s_0N),\,\, v\mapsto \exp(s_0v)$
is a submersion.
\end{step}

By Step 1 $\exp(sN)$ is a dual leaf 
of maximal dimension for  small $s>0$. 
Suppose $s_0<\infty$.
The inclusion $\fL^\#(\exp(s_0x))\subset \exp(s_0N)$ 
follows from $\fL^\#(\exp(sx))\subset \exp(sN)$ for $s<s_0$ 
and Proposition~\ref{prop: general dual}. 

We next want to prove that $\exp(s_0N)$ is contained 
in a dual leaf. Here the definition 
of $\fF^\#$ enters the proof once more. 
 Fix $x\in N$ and let $y\in N$ be any other 
point. 
Choose a piecewise horizontal geodesic 
$\tc$ from the foot-point of $x$ to the foot-point of $y$ 
such that $x$ and $y$ are parallel along $\tc$.
Let $X(t)$ be the parallel vector field along $\tc$ with 
$X(0)=x$.   
By the previous considerations $c_s=\exp(s X(t))$ 
is a variation of curves that maps 
to the trivial variation on a local quotient, $s<s_0$. 
Since $c_0$ projects to a locally minimizing curve 
in a local quotient 
we deduce from the Rauch II comparison theorem and the equality discussion 
that $c_s=\exp(s X(t))$ is a piecewise horizontal geodesic 
as well, $s < s_0$. By continuity the same holds 
for $s=s_0$ and  $\exp(s_0y)$ is contained 
in the same dual leaf as $\exp(s_0x)$.
In other words, $\exp(s_0N)\subset \fL^\#(\exp(s_0x))$.
Thus $\exp(s_0N)=\fL^\#(\exp(s_0x))$. Since $s_0$ 
was chosen maximal it is clear that  
$\fL^\#(\exp(s_0x))$ can not have
 maximal dimension.
It remains to check that the
 map $\psi\colon N\rightarrow \exp(s_0 N)$, $v\mapsto 
\exp(s_0v)$ is a submersion.
Put $\folL^\#_1:=\exp(s_0 N)$. We define a map 
$\varphi\colon N\rightarrow \nu(\folL^\#_1)$ by assigning to $x\in N$ 
the normal vector $y=-\tfrac{d}{ds}_{|s=s_0}\exp(sx)$. 
Clearly $\iota$ is an injective immersion. If we let
 $\pi\colon \nu(\folL^\#_1)\rightarrow \folL^\#_1$ 
denote the natural projection, then $\psi=\pi \circ \varphi$.
Thus it suffices to prove that $\pi_{|\varphi(N)}$
is a submersion. 

Consider vector fields $(X_i)_{i\in I}$  as in Proposition~\ref{prop: general dual}, and let $\gD$ denote the diffeomorphism group generated 
by the flows of these vector fields. In particular the orbits 
of $\gD$ are dual leaves. 
If we identify $\nu_p(\folL^\#_1)$ with $T_pM/T_p\folL^\#_1$, we
get a natural action of $\gD$ on the normal bundle $\nu(\folL^\#_1)$. 
It is clear that $\varphi(N)$ is invariant under this action. 
Since  
$\pi_{|\varphi(N)}$ is equivariant with respect to the $\gD$-action 
it follows that it is  a submersion.

\begin{step} $\folF^\#$ is a singular Riemannian foliation. 
\end{step}

Consider again a dual leaf $\folL^\#_0$ of maximal 
dimension. Notice that the closure  $F$ of the immersed 
submanifold $\folL^\#_0$ in $M$ is contained in the tubular 
neighborhood $B_r(\folL^\#_0)$ from Step~\ref{step1}.
In particular we deduce that $F$ decomposes 
into dual leaves of maximal dimension. 

We claim that the set of points in $M$ for which the dual leaves have
maximal dimension is open and dense. In fact for $q\in M$ 
choose a minimal geodesic $c\colon [0,1]\rightarrow M$ 
 from $F$ to $q$. 
We have seen that $\folL^\#(c(0))$ has maximal dimension 
and clearly  $\dc(0)\in\nu(\folL^\#(c(0)))$. If we let $N$ denote 
the leaf of $\dc(0)$ in $\nu(\folL^\#(c(0)))$, then
 $\exp(sN)$ is a leaf of maximal dimension  for $s\in[0,1)$. 
In fact for each $s\in [0,1)$ 
the map $N\rightarrow \exp(sN)$, $x\mapsto \exp(sv)$
 is injective because otherwise $c$ would not be a minimal 
geodesic from $\folL^\#(c(0))\subset F$ to $q$.

We are now ready to verify that $\folF^\#$ is a singular Riemannian foliation. 
Let $q_0\in M$. It suffices to show that 
 each geodesic emanating perpendicularly to $\folL^\#(q_0)$ 
at $q_0$ stays for a short time perpendicularly to the dual leaves. 
We let
$L_r$ denote the component of $\folL^\#(q_0)\cap B_{5r}(p)$ 
with $q_0\in L_r$ for small $r$.
  Since $\folL(q_0)$ is an immersed submanifold 
it is clear that $L_r$ is Lipschitz continuous in $r\in (0,r_0]$ 
with respect to the Hausdorff distance between subsets of $M$. 
We also may assume that the normal exponential map 
of $L_r$ has an injectivity radius $>3r$.

Clearly we can establish our claim 
by verifying the following statement for some $r>0$:
for any dual leaf $\folL^\#$  
the distance function $d(\cdot,L_r)$ is locally constant 
on $B_r(q_0)\cap \fL^\#$. 
As above we can find a leaf $N\subset \nu(\folL^\#_h)$ 
in the normal bundle of a dual leaf of maximal dimension 
such that $\exp(sN)$ is a dual leaf of maximal 
dimension for all $s\in [0,1)$ and $\fL^\#(q_0)=\exp(N)$. 

We choose an element $u\in N$ with $\exp(u)=q_0$ 
and let for $\delta<<r$, $L_r(\delta)$ 
denote the connected component of
 $\folL^\#(\exp(1-\delta)u)\cap B_{5r}(q_0)$. 
We have seen above that 
the map $N\mapsto \folL(\exp(1-\delta)u), x\mapsto \exp( (1-\delta)x)$
is injective and thus there is a local submersion
$L_r(\delta)\rightarrow L_{r+\delta}$ 
that maps $\exp( (1-\delta)x)$ to $\exp(x)$.
In summary we can say, that the Hausdorff distance between 
$L_r$ and $L_r(\delta)$ is proportional to $\delta$.

Therefore it suffices to check the following holds.
Let $L_1$ denote a component of $\folL_1^\#\cap B_{5r}(q_0)$ 
that intersects $B_r(q_0)$ where $\folL_1^\#$ is a dual leaf 
of maximal dimension. 
Then for any other dual leaf $\folL_2^\#$ 
the function $q\mapsto d(L_1,q)$ 
is locally constant on $\folL_2^\#\cap B_r(q_0)$.

 Fix a point in $q'\in \folL_2\cap B_r(q_0)$. 
We plan to show that  $\folL_2 \cap B_r(q_0) \rightarrow \R$, 
 $q\mapsto d(q,L_1)$ attains a local maximum at $q'$. 
Since $q'$ was arbitrary, this will imply that the function is locally 
constant.
  Choose a vector $v\in TM$ of minimal length 
with a foot point in the closure of $L_1$ 
and with $\exp(v)=q'$. The foot point $p_v$ of 
$v$ is clearly contained in $B_{3r}(q_0)$ and
the dual leaf $\folL^\#(p_v)$
has maximal dimension as well. Furthermore 
an intrinsic open neighborhood $L'$ of $p(v)$ in 
$\folL^\#(p_v)$ is contained in the closure of $L_1$ in 
$M$. In particular $d(q,L_1)\le  d(q,L')$ 
for all $q\in M$. Therefore it suffices 
to prove that $L_2\cap B_r(q_0) \rightarrow \R$, 
 $q\mapsto d(q,L')$ attains a local maximum at $q'$.
Let $N$ be the induced leaf in the normal bundle 
of $\folL^\#(p_v)$ with $v\in N$, and let $N'$ be
 the connected component 
of $N$ intersected with the normal bundle 
of $L'$ with $v\in N'$. 
For all $s<1$ the set $\exp(sN')$ is not contained in a singular
 dual leaf because 
otherwise the geodesic  $\exp(\tau v)$ would not be a minimal connection from 
$L'$ to $q'$. By our previous considerations it follows that for 
all $s\le 1$ the set $\exp(sN')$ is an open subset of a dual leaf.
In particular, $L_2\cap B_r(q_0) \rightarrow \R$, 
 $q\mapsto d(q,L')$ attains a local maximum at $q'$.

\section{Smoothness of the Sharafutdinov Retraction}\label{sec: smooth}
The aim of this section is to prove Corollary~\ref{cor: smooth}.

We consider the Sharafutdinov retraction $P\colon M\rightarrow \Sigma$. 
By Perelman  $P$ is a Riemannian submersion of class $C^{1,1}$. 
Moreover, $P\circ \exp\colon \nu(\Sigma)\rightarrow \Sigma$ 
equals the natural projection $\pi$  from the normal bundle 
 $\nu(\Sigma)$ to the soul $\Sigma$.
 We let $\folF$ denote the fiber decomposition given by $P$ 
and $\folF^\#$ the dual foliation. 
There is a distance tube $B_r(\Sigma)$
 of radius $r$ around $\Sigma$ on which $P$ is of class 
$C^\infty$. Also any horizontal curve 
in $M$ has constant distance to the soul. Thus there is a natural subdivision 
of $B_r(\Sigma)$ into dual leaves. These submanifolds 
are of class $C^\infty$ 
and for suitable small $r$ they are also intrinsically complete 
since they are  via the exponential map diffeomorphic 
to the corresponding dual leaves in $\nu(\Sigma)$.

\begin{thm}\label{lem: commute}\label{lem: vertical geodesics}
Consider the dual foliation $\folF^\#$ 
of an open nonnegatively curved manifold $M$.
Suppose $\folL^\#$ is a dual leaf of class $C^\infty$, and  assume 
that $P_{|\folL^\#}$ is smooth as well. 
\begin{enumerate}
\item[a)] For each $v\in\nu(\folL^\#)$, the curve
 $P(\exp(tv))$ is constant in $t$.
\item[b)]
Let $c(t)\in \folL^\#$
 be a piecewise geodesic which is horizontal with respect 
to $P$, and let $X(t)$ be a parallel vector field along $c$ 
with $X(0)\in\nu(\folL^\#)$.
 Then $\exp(X(t))$ is a piecewise horizontal geodesic with 
respect to $P$ as well. 
\item[c)] 
 Let $F_1=P^{-1}(p_0), F_2=P^{-1}(q_0)$ be fibers of the Sharafutdinov 
retraction. Consider a broken geodesic in $\Sigma$ 
from $p_0$ to $q_0$ and $l\colon F_1\rightarrow F_2, p\mapsto c_p(1)$, 
where $c_p$ denotes the unique horizontal lift of $c$ with $c_p(0)=p$.
For $p\in F_1\cap \folL^\#$ and $q=c_p(1)$ the diagram 
\begin{eqnarray*}
\nu_p\bigl(\folL^\#\bigr) \stackrel{Par_{c_p}}\longrightarrow
 \nu_q\bigl(\folL^\#\bigr)\\
\exp\downarrow\hspace*{5em} \downarrow\exp \\
    F_1\hspace{1em} \stackrel{l}{\longrightarrow}\hspace{1em} F_2\hspace*{2em}
\end{eqnarray*}
commutes, where $\Par_{c_p}$ denotes the parallel transport along $c_p$.
\end{enumerate}
\end{thm}

\begin{sublem} Let $c(t)$ be a horizontal geodesic in 
$M$. Then there is an $(n-1)$-dimensional family $\bbJ$ 
of normal Jacobi fields along $c$ with a self adjoint Riccati operator 
such that each Jacobi field in $\bbJ$ is the variational vector field 
of a variation of horizontal geodesics.
\end{sublem}

If $P$ is of class $C^\infty$ in a neighborhood of $c(0)$, 
 then the sublemma is a general statement 
on Riemannian submersions. Since any horizontal geodesic is a 
limit of such geodesics, the result follows.

\begin{proof}[Proof of Theorem~\ref{lem: commute}.] 
Let $c(t)$ be a horizontal geodesic in $\folL^\#$. 
Choose a family of Jacobi fields $\bbJ$ as in the sublemma. 
As in the proof of Theorem~\ref{thm2} one can show for each 
$J\in \bbJ$ that  if $J(t)\in T\folL^\#$ for some $t$, then 
$J(t)\in  T\folL^\#$ 
 for all 
$t$. As before we deduce
that the normal bundle of $\folL^\#$ along $c$ is spanned 
by parallel Jacobi fields contained in $\bbJ$.

Therefore  each normal vector 
 $v$ of $\folL^\#$ has the property 
that parallel 
transport along any broken $P$-horizontal 
geodesic maps $v$ to a vector which is perpendicular to the dual 
leaf and hence vertical with respect to $P$.
Using Theorem 3.1 in [Guijarro, 2000] part a) and b) follow. 
Part c) is just a simple consequence of b).
\end{proof}

The proof of Theorem 3.1 in  [Guijarro, 2000] is 
a generalization of Perelman's proof of the soul 
conjecture.
 A similar generalization will be given in section~\ref{sec: complete}.

\begin{lem}\label{lem: dual smooth}\label{lem: smooth dual}
The dual leaves are  immersed submanifolds of class $C^\infty$, 
and the restriction of $P$ to each dual 
leaf $\folL^\#_1$ is of class $C^\infty$.
\end{lem}
\begin{proof} As before we choose $r>0$ 
such that $P$ is of class $C^\infty$ in $B_r(\Sigma)$. 
Let $\folL^\#$ a generic dual leaf in $B_r(\Sigma)$. 
In other words the intersection of $\folL^\#$ 
with a fiber of $P$ corresponds to a principal orbit 
of the action of the normal holonomy group on the fiber. 
Then $\folL^\#$ is of class $C^\infty$.
Furthermore the trivialization of the normal 
bundle $\nu(\folL^\#)$ of $\folL^\#$ which is given by Bott 
parallel vector fields is of class $C^\infty$ as well. 
We recall that Bott parallel vectorfield in $\nu(\folL^\#)$
is locally given as the horizontal  lift of a fixed vector in a local quotient 
$U/\folF^\#$ space of the dual foliation.
Consider a Bott parallel vector field $X$ in the normal bundle 
of $\folL^\#$.  Then $X$ is parallel along any horizontal 
geodesic in $\folL^\#$ and by Theorem~\ref{lem: commute}
the image of the map
 \[
h\colon \folL^\#\rightarrow M,\,\, p\mapsto \exp(X(p))
\]
is a dual leaf $\folL^\#_1$.
 Of course it is also clear that all 
dual leaves arise in this way. Moreover the map is of 
class $C^\infty$.
 In order to show that 
$\folL^\#_1$ is of class 
$C^\infty$ it suffices to show 
the map $h$ has constant rank.

The differential of $h$ at $p$ gives rise to a family 
of Jacobi fields along the geodesic $s\to \exp(sX(p))$. 
By Theorem~\ref{lem: commute} this family is the sum of a subfamily of
 parallel Jacobi fields 
which are horizontal with respect to $P$ and a vertical family.
Thus the kernel of $h_{*p}$ is vertical with respect to $P$. 
Let $q$ be another point in $\folL^\#$, $c_p$ be  
horizontal broken geodesic from $p$ to $q$, and put $c= P\circ c_p$. 
Finally we define \[
l\colon P^{-1}(P(p))\rightarrow P^{-1}(P(q))
\]
as in the Theorem~\ref{lem: commute} c).
 Since $P$ is of class $C^{1,1}$ 
the map $l$ is locally bilipschitz. Furthermore $P$ and $l$ 
are of class $C^\infty$ in 
a neighborhood of $p$. By Theorem~\ref{lem: commute} the diagram
\begin{eqnarray*}
 \hspace*{3em}\folL^\#\cap P^{-1}(P(p)) &\stackrel{l}\longrightarrow&
 \folL^\#\cap P^{-1}(P(q))\\
h\downarrow\hspace*{1em}&&\hspace{1em} \downarrow h \\
    F_1\hspace{1em} &\stackrel{l}{\longrightarrow}&\hspace{1em} F_2\hspace*{2em}
\end{eqnarray*}
commutes. Thus
the kernel of $h_{*q}$ is given by the image of the 
kernel of $h_{*p}$ under $l_{*p}$. 
In particular, $h$ is a map of constant rank.
Thus $\folL^\#_1$ is of class $C^\infty$. 
In order to show that $P_{|\folL^\#_1}$ is of class 
$C^\infty$, we observe that $P\circ h=P_{|\folL^\#}$ 
by Theorem~\ref{lem: vertical geodesics}. 
Since $\folL^\#$ is of class $C^\infty$ 
and  
$h\colon \folL^\#\rightarrow \folL^\#_1$ 
is a smooth submersion, it follows that 
 $P_{|\folL^\#_1}$ is of class 
$C^\infty$ as well.
\end{proof}

\begin{proof}[Proof of Corollary~\ref{cor: smooth}]
Let  $p\in M$. By Lemma~\ref{lem: dual smooth}
$\folL^\#(p)$ is of class $C^\infty$ and $P_{|\folL^\#(p)}$ 
is of class $C^\infty$ as well. 
Because of Theorem~\ref{lem: vertical geodesics} 
$P\circ \exp_{\nu(\folL^\#(p))}=P\circ \pi$, where
$\pi\colon \nu(\folL^\#(p))\rightarrow \folL^\#(p)$ is the natural 
projection. Since $P\circ \pi$ is of class $C^\infty$ 
and $\exp_{\nu(\folL^\#(p))}$ is a local diffeomorphism in a neighborhood 
of $0_{p}$, it follows that $P$ is of class $C^\infty$ in a neighborhood 
of $p$.
\end{proof}

\section{Completeness of Dual Leaves.}\label{sec: complete}

This section is devoted to the proof of Theorem~\ref{thm3}.
We first consider the case a). This case is in fact rather obvious.
 Since the group acts transitively on the space of dual leaves, all
 dual leaves have the same dimension 
and, by Proposition~\ref{prop: general},
 the dual foliation is an actual non-singular
 foliation.

{\em b).} Let $\folF$ 
be a Riemannian foliation by $k$-dimensional leaves 
 of a nonnegatively curved
$n$-dimensional compact manifold $M$.
We choose a finite  foliated atlas consisting of 
maps 
$x_i\colon U_i\rightarrow \fD^k\times \fD^{n-k}$ 
where $\fD^k,\fD^{n-k}$ are unit discs in $\R^k$ and 
$\R^{n-k}$ respectively. 

Notice that  for each $i$ the disc $\fD^{n-k}$ 
carries a natural quotient metric $g_i$.
We let $\sigma_i\colon U_i\rightarrow (\fD^{n-k},g_i)$ 
denote the Riemannian submersion.

Choose $\eps>0$ such that  
the injectivity radius of the 
normal exponential map of each of these $k$-dimensional discs 
is larger than $2\eps$ and 
for each point $p\in M$ there is an $i$
with $\bar{B}_{2\eps}(p)\subset U_i$. 
We let $H\subset TM$ denote the set of all unit 
vectors which are perpendicular to the dual 
leaves. By Proposition~\ref{prop: general dual} 
it is clear that $H$ is compact.

We claim that for $v\in H$ the geodesic 
$\exp(sv)$ $(s\in [0,\eps])$ stays in
one leaf of $\folF$.

As in [Guijarro, 2002] we modify Perelman's proof of the soul conjecture
to establish our claim.
We define a displacement function as follows. 
For $v\in H$ consider the foot point $p$ 
of $v$ and  choose an $i$  with $B_{2\eps}(p)\subset U_i$.
Put \[
\dis(s,v):=d\bigl(\sigma_i(p),\sigma_i(\exp(sv))\bigr),
\]
where $\sigma_i\colon U_i\rightarrow (\fD^{n-k},g_i)$ is
 the Riemannian submersion.
It is an important and elementary fact that 
$\dis(s,v)$ is independent of the choice of $i$.
We consider
\[
f(s):=\max\{\dis(s,v)\mid v\in H\}.
\]

Clearly it suffices to prove that  the function 
$f_{[0,\eps]}$ is monotonously decreasing.
Suppose $f(t)>0$ for some  $t\in [0,\eps]$. 
Choose $v\in H$ with $f(t)=\dis(t,v)$ and 
 $i$ with $\bar{B}_{2\eps}(p)\subset U_i$,
 where $p$ is the foot point of $v$.
 \[
f(t)=d(\sigma_i(\exp(tv),\sigma_i(p))).
\]
Let $c\colon [0,1]\rightarrow (\fD^{n-k},g_i)$ 
be the unique minimal geodesic from $\sigma_i(p)$ to $\sigma_i(\exp(tv))$, 
and let $c^h(s)$ be the unique horizontal lift of $c$ starting at $p$.
 By construction there is a $\delta>0$ 
such that the extended geodesics $c_{[-\delta,1]}^h$ and $c_{[-\delta,1]}$ 
are minimal. Furthermore by choosing $\delta$ sufficiently small 
we may assume that $B_{2\eps+\delta f(t)}(p)\subset U_i$. 
Extend $v$ to a parallel vector field $X$ along $c^h$. 
From the proof of Theorem~\ref{thm: micro} we know that 
$X$ stays perpendicular to the dual leaf.

 By applying Rauch comparison
we see that the curve 
$\exp(t X(s))$ ($s\in [-\delta,0]$) is not longer 
than the curve $c_{|[-\delta,0]}$. Thus 
\[
d(\sigma_i(\exp(tX(-\delta))),c(1)))\le 
d(\exp(tX(-\delta)),\exp(tX(0)))\le d(c(0),c(-\delta)).
\]
Therefore 
\begin{eqnarray*}
d(\sigma_i(\exp(tX(-\delta))),c(-\delta))&\ge&
 d(c(1),c(-\delta))-d(c(0),c(-\delta))\\
&=& d(c(1),c(0)).
\end{eqnarray*}
Using $c(-\delta)=\sigma_i(c^h(-\delta))$ our choice of $v$ implies 
 that equality must hold. 
By the equality discussion in Rauch II
 the strip $\exp(\tau X(s))$, $s\in [-\delta,0]$, $\tau\in[0,t]$ is flat. 
Thus
\begin{eqnarray*}
d(\sigma_i(\exp((t-h)X(-\delta))),c(1))^2&\le&
 d(\exp((t-h)X(-\delta)),\exp(tX(0)))^2\\
&\le& d(c(0),c(-\delta))^2+h^2
\end{eqnarray*}
and 
\begin{eqnarray*}
f(t-h)&\ge& d(\sigma_i(\exp((t-h)X(-\delta))),c(-\delta))\\
&\ge& d(c(1),c(-\delta))-d(c(-\delta),c(0)
)- \tfrac{h^2}{2d(c(0),c(-\delta))}\\
&=& f(t)-\tfrac{h^2}{2d(c(0),c(-\delta))}.\\
\end{eqnarray*}
Therefore 
\[
\lim_{h \uparrow 0}\tfrac{f(t)-f(t-h)}{h}\le 0.
\]
Consequently $f_{[0,\eps]}$ is monotonously decreasing and thereby 
constant.

Using the equality discussion in Rauch II 
we see that for a piecewise $\folF$-horizontal
geodesic $c$ in a dual leaf $\folL^\#$
and a parallel unit vector field $X$ along $c$ 
which is normal to $\folL^\#$  the curves 
$t\mapsto \exp(sX(t))$ are piecewise 
$\folF$-horizontal geodesics as well, ($s\in [0,\eps]$).

Suppose there is a dual leaf $\folL^\#$ 
which is not complete. We may assume that $\folL^\#$ 
has minimal dimension among all non-complete leaves.
 Since the intrinsic boundary of $\folL^\#$ 
in $M$ is a union of dual leaves, we can find a dual leaf $\folL^\#_1$  
in the closure of $\folL^\#$ whose dimension with 
$\dim(\folL^\#_1)<\dim(\folL^\#)$.
From the previous claim it is clear that for any $\eps'\le \eps$ 
the $\eps'$ neighborhood around $\folL^\#_1$ is union 
of dual leaves. 
By construction $\folL^\#_1$ is in the closure 
of $\folL^\#$. Since $\folL^\#$ and $\folL^\#_1$ 
have different dimensions, we can
 employ Proposition~\ref{prop: general dual} to see that
 $\folL^\#$ and $\folL^\#_1$ have positive Hausdorff 
distance in $M$. Combining the last three statements gives
a contradiction.

For the proof of c) notice that we can apply
 Theorem~\ref{lem: vertical geodesics} to 
see that  the $\eps$-neighborhood 
of a complete dual leaf decomposes into dual leaves for all small $\eps>0$. 
 As in the 
previous paragraph this gives the completeness of all dual 
leaves.

\section{Totally Geodesic Flats in Foliated  
Manifolds with Nonnegative Sectional Curvature.}\label{sec: flats}

In this section we prove
 Corollary~\ref{cor: flats}. In fact it
 clearly follows from the following more 
general result.
\begin{prop}\label{prop: flats}
 Let $\folF$ be a singular Riemannian foliation of 
a nonnegatively curved manifold $M$ and 
suppose the dual foliation $\folF^\#$ has complete 
leaves. Let $x\in T_pM$ be a vector that 
is horizontal with respect to $\folF$ and 
$v\in T_pM$ a vector that is horizontal with 
respect to $\folF^\#$. Then $x$ and $v$ span a totally geodesic flat.
\end{prop}

\begin{proof} Let $c(t)=\exp(tx)$ and let 
$V(t)$ be parallel along $c$ with $V(0)=v$. 
We have seen that $V(t)$ stays perpendicular to the dual 
leaf $\folL^\#$. Furthermore for 
each $t$ the curve $s\mapsto \exp(sV(t))$ is a horizontal geodesic 
with respect to $\folF^\#$ and hence it is vertical 
with respect to $\folF$. 
For each $t$ one can find an $\eps>0$ 
such that $c_{|[t,t+\eps]}$ is a local minimal connection 
between $\folL(c(t))$ and $\folL(c(t+\eps))$. 
By Rauch II the parallel curves 
 $t\mapsto \exp(sV(t))$ are not longer. 
Since these curves connect the same leaves 
equality must hold in Rauch's comparison theorem
 and thus $c$ and $V$ generate 
a totally geodesic flat.
\end{proof}

\section{Non-contractible, Nonnegatively Curved Open  Manifolds 
have nontrivial Products as Metric Quotients.}\label{sec: product}

In this section we prove
Corollary~\ref{cor: product submetry}.

\begin{prop}\label{prop: intersection}
 Let $\folF$ be a singular Riemannian foliation of 
a nonnegatively curved manifold $M$ and 
suppose the dual foliation $\folF^\#$ has complete 
leaves. We define a singular foliation 
$\folF\cap \folF^\#$ by the property that the leaf of a point $p$ 
is given by the $p$-component of $\folL(p)\cap \folL^\#(p)$. 
Then $\folF\cap \folF^\#$ is a transnormal system.
\end{prop}

\begin{proof}
 Notice that $\folL(p)$ and $ \folL^\#(p)$
intersect transversely. So $\folF\cap \folF^\#$ 
is indeed a subdivision into intrinsically complete 
immersed submanifolds. Let $u\in T_pM$ be perpendicular to 
$ \folL(p)\cap \folL^\#(p)$. 
Then $u=x+v$ with $x\in \nu_p(\folL(p))$ and $v\in \nu_p(\folL^\#(p))$. 
By Proposition~\ref{prop: flats}
 $x$ and $v$ span a totally geodesic flat. 
Moreover it is clear form the proof of Proposition~\ref{prop: flats} 
that at each point the flat is spanned by one $\fF$-horizontal 
 and one $\fF^\#$-horizontal vector. 
Therefore all tangent vectors of the flat are $\fF\cap \fF^\#$-horizontal
and hence the same holds for the curve $\exp(tu)$.
\end{proof}

\begin{proof}[Proof of Corollary~\ref{cor: product submetry}.]  We let $\folF$ denote the foliation induced 
by the Sharafutdinov retraction and $\folF^\#$ its dual.
We define the leaves 
of $\overline{\folF}^\#$  as the closures of  leaves 
of $\folF^\#$. Clearly the leaves of $\overline{\folF}^\#$
are the fibers of a globally defined proper submetry
$\sigma_2\colon M\rightarrow A$, where $A$ is
 a noncompact Alexandrov space.

Furthermore there is a distance tube $B_r(\Sigma)$ around the soul 
such that the leaves are via the exponential map isomorphic to 
the corresponding closures of dual 
leaves in $\nu(\Sigma)$. In particular the leaves of 
 $\overline{\folF}^\#$ in $B_r(\Sigma)$ are of class $C^\infty$. 
Analogously to the proof of Lemma~\ref{lem: smooth dual} 
one can now show that all leaves in  $\overline{\folF}^\#$
are of class $C^{\infty}$.

Thus $\overline{\folF}^\#$ is 
a transnormal system. 
As in the proof of Proposition~\ref{prop: intersection} 
one can show that $\overline{\folF}^\#$ intersected 
with the fibers of $P$ gives a transnormal system as well. 
Hence the map 
$\sigma:=(P,\sigma_2)\colon M\rightarrow \Sigma\times A$ 
is a submetry. Clearly the fibers of $\sigma$ 
are compact and smooth.
\end{proof}

 Notice that the fibers of $\sigma$
 are given by the closures of 
orbits  
of the normal holonomy group of the soul acting on the fibers of $P$.

\section{The Horizontal Distribution of an Open 
Nonnegatively Curved Manifold is Linear.}\label{sec: linear}

In this section we prove Corollary~\ref{cor: linearization}.
We start with an observation that is somewhat related to the construction
 of the Sharafutdinov retraction.

\begin{lem}\label{lem: sharafut}
 Let $M$ be an open nonnegatively curved 
manifold and let $\folF^\#$ denote the dual foliation of the 
Sharafutdinov retraction $P\colon M\rightarrow \Sigma$. 
\begin{enumerate}
\item[a)] The convex exhaustion obtained from the soul construction 
is invariant under the dual foliation.
\item[b)] For each dual leaf $\folL^\#\neq \Sigma$, 
there is a sequence of dual leaves 
$\folL^\#_n$ converging to $\folL^\#$ with 
$\dim(\folL^\#_n)=\dim(\folL^\#)$,
 $\dim(\bar{\folL}^\#_n)=\dim(\bar{\folL}^\#)$
and $d(\Sigma,\folL_n^\#)<d(\Sigma,\folL^\#)$.
\end{enumerate}
\end{lem}
\begin{proof} 
{\em a).}
We start by considering the Busemann function of a point 
$p_0\in M$, $b(x)=\lim_{r\to \infty} d(\partial B_r(p_0),x)-r$.

Let $c\colon \R\rightarrow M$ be a horizontal geodesic. 
Then $c$ is contained in a relatively compact set and thus 
$b\circ c$ is bounded. On the other hand $b\circ c$ is concave and 
hence $b\circ c$ is constant.
This simple observation shows that the levels of $b$ decompose into
 dual leaves.

Let $C$ denote the maximal level of the Busemann function, and let $\partial 
C$ denote the intrinsic boundary of $C$.
Let $\Sigma$ denote the soul of $C$.
If $\partial C$ is empty then $C=\Sigma$ and we are done.

We have just seen that $C$ decomposes into dual leaves. 
Notice that the dual leaf $\Sigma$ has constant distance to $\partial C$. 
Put
\[ 
G:=\{p\in C- \partial C\mid \folL^\#(p) \mbox{ has constant 
distance to $\partial C$}\}.
\]
We have just seen that $G$ is not empty and clearly $G$ is closed in 
$C-\partial C$. We claim that $G$ is open 
in $C$ as well. Let $\folL^\#$ denote a dual leaf in $G$ 
and let $r$ denote the distance to $\partial C$. 
Then the set $B_r(\folL^\#)\cap C$ decomposes into dual leaves. 
Let $c(t)$ be a horizontal geodesic in $B_r(\folL^\#)\cap C$. 
 As before it is clear that 
$t\mapsto d(\partial C,c(t))$ is both bounded and concave 
and thereby 
constant. 
Thus all dual leaves in  $B_r(\folL^\#)\cap C$ have constant distance 
to $\partial C$ and this in turn shows that $G$ is open in $C$.

We have proved 
 that the level sets of $d(\partial C,\cdot)$ 
decompose into dual leaves as well. A simple induction 
argument shows that the whole convex exhaustion is invariant 
under the dual foliation.

{\em b).}
For each point $p\in M\setminus \Sigma$ there is a unique 
convex set $C$ in the convex exhaustion such that $p$ is contained 
in the intrinsic boundary $\partial C$ of $C$. 
 By a) $\folL^\#(p)\subset \partial C$. 
We let $T_pC$ denote the tangent cone of $C$, and
let $c$ 
be a horizontal geodesic in $\folL^\#(p)$.
Since $C$ decomposes into dual leaves we can employ 
Theorem~\ref{lem: commute} to see that $T_{c(t)}C\cap\nu_{c(t)}(\folL^\#(p))$ 
is parallel along $c$.
For each $q\in \folL^\#(p)$ we define $X_q$ as the unique 
unit in the tangent cone $T_qC$ with maximal distance to 
the boundary of $T_qC$. 
Clearly $X$ is normal to the dual leaf and 
 parallel along any horizontal geodesic in $\folL^\#(p)$. 
This proves that for each $s>0$, the image of 
$\folL^\#(p)\rightarrow M, q\mapsto \exp(s X_{|q})$ is a dual leaf 
as well. Clearly its distance to the soul is smaller than the distance 
of $\folL^\#(p)$. Moreover its dimension  constant in
$s$ for small $s$.
\end{proof}

We recall that the normal holonomy group of the soul 
does not need to be compact even in the simply connected case. 

\begin{prop}\label{prop: closure} Let $M$ be an open manifold of nonnegative 
sectional curvature, $\Sigma$ a soul of $M$ and $p\in \Sigma$.
 Consider a fiber $F:=P^{-1}(p)=\exp(\nu_p(\Sigma))$
 of the Sharafutdinov retraction. 
The normal holonomy group $\gH$ of $\Sigma$ 
 acts on $F$ by diffeomorphisms and
the image of the 
induced homomorphism $\gH\rightarrow \Diff(F)$ 
has a compact Lie group as its closure.
\end{prop}

\begin{proof}

Let $F_r$ denote the ball of radius 
$r$ in $F$ around $p$. 
Notice that the action of $\gH$ leaves  $F_r$ 
invariant. 
For small $r$ it is clear that 
the homomorphism $\gH\rightarrow \Diff(F_r)$ 
has a relatively compact image, since 
the action of $\gH$ is via the exponential map 
isomorphic to a linear orthogonal action. 

Choose $r\in (0,\infty]$ maximal such that the image 
of the above homomorphism is relatively compact. 
Suppose, on the contrary, that $r<\infty$.

Even though the action of $\gH$ is not isometric, we can use 
Theorem~\ref{lem: commute} to see that
$\|L_{h*} v\|=\|v\|$ and 
$h \exp(v)=\exp(L_{h*}v)$ for any vector $v\in \nu_q(Hq)$ 
in the normal bundle of an orbit, where $L_h(q):=hq$. 
Consider an orbit $\bgH\star q$ of the closure $\bgH\subset\Diff(F_r)$. 
Let $\eps$ denote the focal radius of the normal exponential map 
of $\bgH\star q$. The above discussion shows 
that the $\gH$ action in the tubular neighborhood 
$B_{\eps}(\bgH \star q)=B_{\eps}(\gH\star q)$ 
extends naturally to an action of $\bgH$.
Thus it suffices to show that the union of 
all tubes $B_{\eps}(\gH\star q)$ covers the closure of $F_r$.

By Lemma~\ref{lem: sharafut} each closure of an $\gH$ orbit
 in the boundary of $F_r$ can be approximated by a sequence of closures 
of $\gH$ orbits in $F_r$ 
which have the same dimension as the given one. 
Since the closure of these orbits are the smooth fibers of a submetry  on
$F$, it follows that the focal radii of these submanifolds  
stay bounded below 
and thus each orbit in the boundary of $F_r$ 
is contained in some $B_{\eps}(\gH\star q)$ with $q\in F_r$.
\end{proof}

\begin{proof}[Proof of Corollary~\ref{cor: linearization}.]
Consider a fiber $F=P^{-1}(p)$ of the Sharafutdinov retraction. 
Recall  that the distance function of $p$ 
has no critical points in $F$. 
Thus we can find a gradient like unit vector field $X$  
in $F\setminus p$. 

By Proposition~\ref{prop: closure} the 
closure $\bgH$ acts on $F$. For any vector $v\in \nu_q(\bgH \star q)$ 
and any $h\in \bgH$ we have $\| L_{h*}v\|=\|v\|$ and 
$h\exp(v)=\exp(L_{h*}v)$, see Theorem~\ref{lem: commute}. 
Using that the orbits of $\bgH$ induce a singular Riemannian 
foliation we see furthermore that 
 all minimal geodesics 
from $q$ to $p\in \Sigma$ are perpendicular to $\bgH\star q$. 

It is now easy to see that for any $h\in \bgH$ the vector field 
$\tilde{X}_{|q}:=L_{h*}X_{|h^{-1}q}$ is a again a gradient like vector field.
A simple averaging argument shows that we can 
find a gradient like vector field $Y$ of bounded length 
that commutes with the action of $\bgH$. 
We can also assume $Y$ coincides in a small pointed neighborhood 
$B_{\delta}(p)\setminus p$
of $p$ with the actual gradient of the distance function. 
 
Since $Y$ commutes with the action of the holonomy group, 
there is a unique way to extend $Y$
to a vertical gradient like vector field $Z$ on $M$, by pushing 
$Y$ with diffeomorphism as in Theorem~\ref{lem: commute} to different fibers. 
Notice $Z$ is given by Jacobi fields along horizontal 
geodesics and the flow of $Z$ maps 
horizontal geodesics to horizontal geodesics.

We now consider the diffeomorphism 
$f\colon \nu(\Sigma)\rightarrow M$ given as follows: 
for $r v\in \nu(\Sigma)$ with $\|v\|=1$ and $r\ge 0$ 
consider the integral curve $\gamma$ of $Z$ with $\gamma(0)=\exp(\delta v)$ 
and put $f(rv)=\gamma(r-\delta)$. Notice that $f(rv)=\exp(rv)$ for
 $r\le \delta$. Since the vectorfield $Z$ is vertical, $f$ 
satisfies a). Since the flow of $Z$ maps horizontal geodesics 
to horizontal geodesics, $f$ maps parallel vectorfields 
along geodesics in $\Sigma$ onto horizontal geodesics in $M$, 
as claimed in b).
\end{proof}

\begin{rem}\begin{enumerate}
\item[a)] Of course Corollary~\ref{cor: linearization} 
implies that the Alexandrov space $A$ from 
Corollary~\ref{cor: product submetry} is bilipschitz equivalent to
$\nu_p(\Sigma)/\bar{\gH}$, where $\bar{\gH}$ denotes 
the closure of the normal holonomy group of the soul, there is of course
no global bilipschitz constant though.
\item[b)]
Because of Lemma~\ref{lem: sharafut} the convex 
exhaustion obtained from the soul construction in $M$ 
is just the inverse 
image of the convex exhaustion of the soul construction in $A$ 
under the submetry $\sigma_2\colon M\rightarrow A$.
\end{enumerate}
\end{rem}

\section{Rigidity of Non-Primitive Actions
 in Nonnegative Sectional Curvature}\label{sec: actions}

This section is devoted to the proof of Corollary~\ref{cor: actions}.

\begin{prop}\label{prop: actions} Suppose a Lie group $\G$ acts isometrically on 
a nonnegatively curved manifold $M$. Let $\folF$ denote the singular 
Riemannian foliation 
induced by the orbit decomposition of $\G$.
 Suppose the dual foliation $\folF^\#$ has a leaf which is not dense. 
Then there is a closed subgroup $\gK\subsetneq \G$,
 an invariant metric $\G/\gK$ 
and a $\G$-equivariant 
Riemannian submersion $\sigma\colon M\rightarrow \G/\gK$.
 Furthermore the fibers of $\sigma$ are closures of leaves 
of $\folF^\#$.
\end{prop}

\begin{proof} Let $\bar{\folF}^\#$ denote the foliation 
whose leaves are given by the closures of leaves in $\folF^\#$.
Clearly the group $\G$ acts transitively on the space of dual leaves 
and hence also on the space of leaves of $\bar{\folF}^\#$. 
Since $\folF^\#$ is a singular Riemannian submersion by Theorem~\ref{thm: micro} and Theorem~\ref{thm: complete}
the leaves of $\bar{\folF}^\#$  are the fibers 
of a submetry $\sigma\colon M\rightarrow X$.  
The action of $\G$ on $M$ induces a transitive isometric action 
of $\G$ on $X$ and hence $X$ is a homogeneous space $\G/\gK$.

Consider the fiber $F:=\sigma^{-1}(\gK)$. 
Then  $F=\exp(\nu(\gK p)\cap \nu(\G \star p))$ for each $p\in F$ 
and hence $F$ is smooth.
Thus $\sigma$ is a Riemannian submersion.
\end{proof}

\begin{proof}[Proof of Corollary~\ref{cor: actions}.]
Consider a fixed point component $N$ of the principal isotropy 
group $\gH$ which intersects a principal orbit.
Using that $M/\G$ has no boundary we deduce from Proposition 11.3 
in [Wilking, 2003] 
that any isotropy group of a point $p\in N$ is contained in 
the normalizer $N(\gH)$ of $\gH$. In particular it follows that 
for all $p\in N$ the normal space $\nu_p(\G \star p)$ 
of the orbit $\G \star p$ is contained in $T_pN$. 
Therefore for each $p\in N$ the dual leaf $\folL^\#(p)$ 
of $p$ is contained in 
$N$. 
By Proposition~\ref{prop: actions} there is a Riemannian 
submersion $\sigma\colon M\rightarrow \G/\gK$. 
Furthermore one fiber $F$ of $\sigma$ is contained in $N$. 
This in turn shows $\gK\subset N(\gH)$.

We next plan to show each fiber $F:=\sigma^{-1}(x)$ 
is totally geodesic. 
For that we consider a dense dual leaf $\folL^\#\subset F$.
Let $M'$ denote the union of all principal orbits in $M$, and  
let $p\in M'\cap \folL^\#$.
Using that $M/\G$ has no boundary it is not hard to check
\[
M'\cap \folL^\#:=\bigl\{q\in M\mid \mbox{ there is piecewise 
horizontal curve in $M'$ from $p$ to $q$}\bigr\}. 
\]
Consider a Killing field $X$ which is perpendicular to $\folL^\#$ 
at $p$. We claim that $X$ is perpendicular to $\folL^\#$ for all 
$q\in \folL^\#$. To prove this we may assume that $q\in \folL^\#\cap M'$. 
Since there is a piecewise horizontal geodesic in $M'$ 
from $p$ to $q$ it suffices to prove the following.
If $c(t)$  ($t\in [0,1]$) is a horizontal geodesic 
in $M'\cap \folL^\#$ and $X_{|c(0)}\in\nu(\folL^\#)$, then  
$X_{|c(1)}\in\nu(\folL^\#)$. 
But this is clear since by the proof of Theorem~\ref{thm: micro}
$X$ is a parallel Jacobi field along $c$ which is perpendicular 
to $\folL^\#$.

This shows that each  Bott parallel vector field $X$ along $\folL^\#$ 
is the restriction of a  Killing field. 
By restricting attention 
to those Bott parallel vector fields along $\folL^\#$ which 
are perpendicular to the closure $F$ of $\folL^\#$, we see that
 the Bott parallel vector fields 
along $F$ with respect to $\sigma$ 
are given by Killing fields. 
Of course this implies that the holonomy maps of the submersion 
$\sigma\colon M\rightarrow \G/\gK$ are isometries and hence the fibers 
of $\sigma$ are totally geodesic.
\end{proof}

\begin{rem} Notice that the fibers of the Riemannian submersion $\sigma$ 
of Proposition~\ref{prop: actions} are pairwise isometric.  
It would be interesting to know whether the fibers 
have nonnegative curvature.
\end{rem}

\section{A Slice Theorem for Dual Foliations} 
In nonnegative curvature dual foliations have an additional
 remarkable property.

\begin{thm} 
Let $\folF$ be a singular Riemannian foliation 
of a nonnegatively curved manifold $M$. 
Suppose the dual foliation has closed leaves. 
Then there there is subgroup $\gD\subset \Diff(M)$ such that 
the leaves of $\folF^\#$ are orbits of the $\gD$--action 
and for each dual leaf $\folL^\#$  the $\gD$ action on a
suitable 
 tubular neighborhood 
$B_r(\folL^\#)$  is orbit equivalent 
to the natural action of $\gD$ on $\nu(\folL^\#)$. 
\end{thm}

We recall that the action of $\gD$ on $\nu(\folL^\#)$ is induced by the 
identification $\nu_q(\folL^\#)=T_qM/T_q\folL^\#$ for 
all $q\in \folL^\#$.
In particular it follows that each tangent cone
of the orbit space $M/\gD$ is isometric 
to $\R^d/\gH$, where $\gH$ is a suitable subgroup of $\Or(d)$.
In other words, the singularities of $M/\gD$ look like 
singularities on an orbit space of an isometric group action.

\begin{proof} Define $\gD$ as in the proof
 of Proposition~\ref{prop: general dual}.
Consider a unit normal vector $v\in \nu_p(\folL^\#)$,
 a piecewise horizontal geodesic in $\folL^\#$
starting at $p$, and the parallel vector field $X$ 
along $c$ with $X(0)=v$. By Proposition~\ref{prop: flats}
and Proposition~\ref{prop: intersection}
 the curve $t\mapsto \exp(s(X(t))$ is a 
piecewise horizontal geodesic as well for each $s\in\R$.

Let $N$ 
denote the subset of $\nu(N)$ consisting of all unit vectors 
which are parallel to $v$ along some piecewise horizontal 
geodesic. 
It follows that for all $s$ the set $\exp(sN)$ is a dual leaf. 
Hence it suffices to prove that $N$ is an orbit of $\gD$ with respect 
to the induced action of $\gD$ on $\nu(\folL^\#)$.

Consider a horizontal geodesic $c$ in $\folL^\#$ 
and a parallel vector field $H$ normal to $\folL^\#$. 
We extend $\dc(0)$ to a vector field $X$ 
along $\folL:=\folL(c(0))$ in neighborhood of $c(0)$ 
by using radially normal parallel translation.
We now choose a vector field $Y$ in the normal bundle $\nu(\folL)$ 
with compact support such that $Y$ contained in the kernel 
of $\pi_*$ and
$Y_{sX_p}=\tfrac{d}{dt}_{t=s} tX$ in a neighborhood of $(p,s)=(c(0),0)$.
We may assume that there is an open set $U\subset \nu(\folL)$ 
containing the support of $Y$ such that $\exp_{U}$ is a diffeomorphism onto 
its image. 

Let $Z$ be the vector field in $M$ which is $\exp$-related to $Y$. 
By construction the flow of $Z$ is contained in $\gD$. 
The induced action of the flow of $Z$ in $\nu(\folL^\#)$ has 
$H_{[-\eps,\eps]}$ as an integral curve for a suitable small $\eps>0$.
In summary it follows that there is an $\eps>0$ such that 
$H([-\eps,\eps])$ is contained in a $\gD$-orbit in $\nu(\folL^\#)$.
A simple compactness argument shows that $N$ is 
 a $\gD$-orbit in $\nu(\folL^\#)$.
\end{proof}

\section*{Final Remarks} 

\begin{rems*}\begin{enumerate}
\item[a)]  One can show that for a transnormal system 
there are Lipschitz continuous vectorfields $\{X_i\mid i\in I\}$ 
such that $T_p\folL(p)=\spann_{\R}\{X_{i|p}\mid i\in I\}$. 
Using this it is clear that Theorems 1 and 2 remain valid 
if one just assumes that $\fF$ is a transnormal system. 
\item[b)] If $\fF$ is a 
 transnormal system such that all leaves have the same dimension, 
then $\fF$ is a (non-singular) Riemannian foliation. 
\item[c)] To the best of the authors knowledge it is not known 
whether there is any transnormal system which is not
a singular Riemannian foliation. There are claims 
in the literature that examples exists, but these 
claims also would imply that part b) of this remark is false. 
\end{enumerate} 
\end{rems*}

The proofs of these remarks are elementary
 but not trivial. Maybe the details will be carried out
 somewhere else.


\hspace*{1em}\\
\begin{footnotesize}
\hspace*{0.3em}{\sc University of M\"unster,
Einsteinstrasse 62, 48149 M\"unster, Germany}\\
\hspace*{0.3em}{\em E-mail address: }{\sf wilking@math.uni-muenster.de}
\end{footnotesize}
\end{document}